\documentclass[final,narrowdisplay]{elsart1p}

\usepackage{graphicx,subfigure,epsfig}
\usepackage{amssymb,latexsym} 
\usepackage{threeparttable,algorithm,algorithmic,caption}

\newcommand{\be}{\begin{eqnarray}}
\newcommand{\ee}{\end{eqnarray}}

\newcommand{\bitem}{\begin{itemize}}
\newcommand{\eitem}{\end{itemize}}

\newcommand{\bE}{\mathbb{E}}

\newcommand{\bZ}{\mathbb{Z}}

\newcommand{\fv}{\mathbf{v}}

\newcommand{\cP}{\mathcal{P}}

\begin{document}
\begin{frontmatter}
\title{Fast Poisson Noise Removal by Biorthogonal Haar Domain Hypothesis Testing} 
\author[A1]{B. Zhang\corauthref{cor}}
\corauth[cor]{Corresponding author.}
\address[A1]{Quantitative Image Analysis Group URA CNRS 2582 of Institut Pasteur, 75724 Paris, France}
\ead{bzhang@pasteur.fr}

\author[A2]{M. J. Fadili} 
\address[A2]{Image Processing Group GREYC CNRS UMR 6072, 14050 Caen Cedex, France}
\ead{Jalal.Fadili@greyc.ensicaen.fr}

\author[A3]{J.-L. Starck}  
\address[A3]{DAPNIA/SEDI-SAP, Service d'Astrophysique, CEA-Saclay, 91191 Gif sur Yvette, France}
\ead{jstarck@cea.fr}

\author[A4]{S. W. Digel} 
\address[A4]{Stanford Linear Accelerator Center, 2575 Sand Hill Road, Menlo Park, CA 94025}
\ead{digel@slac.stanford.edu}

\begin{abstract}
Methods based on hypothesis tests (HTs) in the Haar domain are widely used to denoise Poisson count data. Facing large datasets or real-time applications, Haar-based denoisers have to use the decimated transform to meet limited-memory or computation-time constraints. Unfortunately, for regular underlying intensities, decimation yields discontinuous estimates and strong ``staircase'' artifacts. In this paper, we propose to combine the HT framework with the decimated biorthogonal Haar (Bi-Haar) transform instead of the classical Haar. The Bi-Haar filter bank is normalized such that the $p$-values of Bi-Haar coefficients ($p_{BH}$) provide good approximation to those of Haar ($p_H$) for high-intensity settings or large scales; for low-intensity settings and small scales, we show that $p_{BH}$ are essentially upper-bounded by $p_H$. Thus, we may apply the Haar-based HTs to Bi-Haar coefficients to control a prefixed false positive rate. By doing so, we benefit from the regular Bi-Haar filter bank to gain a smooth estimate while always maintaining a low computational complexity. A Fisher-approximation-based threshold implementing the HTs is also established. The efficiency of this method is illustrated on an example of hyperspectral-source-flux estimation.
\end{abstract}
\end{frontmatter}

\section{Introduction}
\label{sec:intro}
Astronomical data analysis often requires Poisson noise removal \cite{Starck1998book}. This problem can be formulated as follows: 
we observe a $q$-dimensional ($q$D) discrete dataset of counts $\fv=(v_i)_{i\in \bZ^q}$ where $v_i$ follows a Poisson distribution of intensity $\lambda_i$, i.e.~$v_i \sim \cP(\lambda_i)$. Here we suppose that $v_i$'s are mutually independent. The denoising aims at estimating the underlying intensity profile $\Lambda=(\lambda_i)_{i\in \bZ^q}$ from $\fv$.

A host of estimation methods have been proposed in the literature (see the reviews \cite{Besbeas2004}\cite{Willett06} and their citations), among which an important family of approaches based on hypothesis tests (HTs) is widely used in astronomy \cite{Kolaczyk1997,Kolaczyk2000}\cite{Charles2003}. These methods rely on Haar transform and the HTs are applied on the Haar coefficients to control a user-specified false positive rate (FPR). When working with large datasets or real-time applications, the decimated Haar transform is generally required to meet limited-memory or computation-time constraints. This is even more true when processing astronomical hyperspectral data, which are usually very large in practice. Unfortunately, for regular underlying intensities, decimation yields discontinuous estimates with strong ``staircase'' artifacts, thus significantly degrading the denoising performance. Although \cite{Bijaoui2001} and \cite{Kolaczyk1999a} attempted to generalize the HTs for wavelets other than Haar, \cite{Bijaoui2001} is more computationally complex than Haar-based methods, and \cite{Kolaczyk1999a} adopts an asymptotic approximation which may not allow reasonable solutions in low-count situations. In an astronomical image decompression context, \cite{Bobichon1997} has also proposed to remove Haar block artifacts by minimizing at each resolution level the $\ell^2$-norm of the Laplacian of the solution under some constraints on its wavelet coefficients. It has been shown that this approach was efficient in removing the artifacts, but it requires solving $J$ minimization problems, where $J$ is the number of scales. This can be quite time-consuming and would limit the interest in using Haar for large-dataset analysis.

In this paper, we propose to combine the HT framework with the decimated bi-orthogonal Haar (Bi-Haar) transform. The Bi-Haar filter bank is normalized such that the $p$-values of Bi-Haar coefficients ($p_{BH}$) approximate those of Haar ($p_H$) for high-intensity settings or large scales; for low-intensity settings and small scales, we show that $p_{BH}$ are essentially upper-bounded by $p_H$. Thus, we may apply the Haar-based HTs to Bi-Haar coefficients to control a prefixed FPR. By doing so, we benefit from the regular Bi-Haar filter bank to gain a smooth estimate. A Fisher-approximation-based threshold implementing the HTs is also established. We find that this approach even exhibits a performance comparable to the more time/space-consuming translation-invariant Haar (TI Haar or undecimated Haar) denoising in some of our experiments. The efficiency of this method is also illustrated on an example of hyperspectral-source-flux estimation.

The paper is organized as follows. We begin with the review of the wavelet HTs in Section \ref{sec:hypwav}, and then Bi-Haar domain tests are presented in Section \ref{subsec:bihaar}. Section \ref{subsec:thresh} details some thresholding operators implementing the tests. The final denoising algorithm is summarized in Section \ref{subsec:algo}, and the numerical results are shown in Section \ref{sec:res}. We conclude in Section \ref{sec:conc}, and the mathematical details are deferred to the appendices.

\section{Hypothesis testing in the wavelet domain\label{sec:hypwav}}
Wavelet domain denoising can be achieved by zeroing insignificant coefficients while preserving significant ones. We detect significant coefficients by applying a binary HT on each wavelet coefficient $d$:
\[ H_0: d=0\ \mbox{vs.}\ H_1: d\neq 0 \]
Note that since any wavelet has a zero mean, if $d$ comes from a signal of constant intensity within the wavelet support, then $d \in H_0$.

Individual HTs are commonly used to control a user pre-specified FPR in the wavelet domain, say $\alpha$. The tests are carried out in a coefficient-by-coefficient manner. That is, the $p$-value of each coefficient $p_i$ is calculated under the null hypothesis $H_0$. Then, all the coefficients with $p_i>\alpha$ will be zeroed. If we desire to control global statistical error rates, multiple HTs may be adopted such as Bonferroni correction which controls the Family-Wise Error Rate (FWER), and the Benjamini and Hochberg procedure \cite{Benjamini1995}\cite{Benjamini2001} controlling the false discovery rate (FDR).

\subsection{$p$-values of wavelet coefficients under $H_0$\label{subsec:pvalue}}
To carry out HTs, we need to compute the $p$-value of each wavelet coefficient under $H_0$. Although the probability density function (\textit{pdf}) of a $H_0$-coefficient has been derived in \cite{Bijaoui2001}, this \textit{pdf} has no closed form for a general wavelet. Thus the $p$-value evaluation in practice is computationally complex. 

To obtain distributions of manageable forms, simple wavelets are preferred, such as Haar. To the best of our knowledge, Haar is the only wavelet yielding a closed-form \textit{pdf}, which is given by \cite{Skellman1946} ($n\geq 0$): 
$\Pr(d=n; \lambda) = e^{-2\lambda}I_n(2\lambda)$, where $d=X_1-X_2$, $X_1, X_2\sim \mathcal{P}(\lambda)$, and $I_n$ is the $n$-th order modified Bessel function of the first kind. For negative $n$, the probability can be obtained by symmetry. The tail probability ($p$-value) is given by \cite{Johnson1959}:
\begin{equation}
\label{eq:haartail}
\Pr(d \geq n; \lambda) = \Pr\left(\chi_{(2n)}^2(2\lambda) < 2\lambda\right), \quad n \geq 1 
\end{equation}
where $\chi^2_{(f)}(\Delta)$ is the non-central chi-square
distribution with $f$ degrees of freedom and $\Delta$ as non-centrality parameter. 

\subsection{Bi-Haar domain testing\label{subsec:bihaar}}
Haar wavelet provides us with a manageable distribution under $H_0$. But due to the lack of continuity of Haar filters, its estimate can be highly irregular with strong ``staircase'' artifacts when decimation is involved. 

To solve this dilemma between distribution manageability and reconstruction regularity, we propose to use the Bi-Haar wavelet. Its implementation filter bank is given by \cite{Starck1998book}:
\[
\begin{array}{ll} 
h = 2^{-c}[1,1], \quad & g = 2^{-c}r[\frac{1}{8},\frac{1}{8},-1,1,-\frac{1}{8},-\frac{1}{8}];\\
\tilde{h} = 2^{c-1}r[-\frac{1}{8},\frac{1}{8},1,1,\frac{1}{8},-\frac{1}{8}], \quad & \tilde{g} = 2^{c-1}[1,-1]
\end{array} 
\] 
where $c$ and $r = (1+2^{-5})^{-1/2}$ are normalizing factors, $(h,g)$ and $(\tilde{h},\tilde{g})$ are respectively the analysis and synthesis filter banks. Note that our Bi-Haar filter bank has an unusual normalization. The motivation behind this is to ensure that the Bi-Haar coefficients will have the same variance as the Haar ones at each scale. Let us also point out that to correct for the introduction of the factor $r$, the Bi-Haar coefficients must be multiplied by $r^{-1}$ at each stage of the recursive reconstruction. For comparison, the Haar filter bank is ($h = 2^{-c}[1,1]$, $g = 2^{-c}[-1,1]$, $\tilde{h} = 2^{c-1}[1,1]$, $\tilde{g}=2^{c-1}[1,-1]$). It follows that the synthesis Haar scaling function is discontinuous while that of Bi-Haar is almost Lipschitz~\cite{Villasenor1995}\cite{Rioul1992}. Hence, the Bi-Haar reconstruction will be smoother.

At scale $j\geq 1$, let us define $\lambda_j=2^{j}\lambda$ where $\lambda$ is the underlying constant intensity. Then, a Haar coefficient can be written as $d^h_{j} = 2^{-cj}(X_1-X_2)$ where $X_1,X_2\sim\cP(\lambda_j/2)$ are independent. We note $p_H := \Pr(d^h_j \geq 2^{-cj}k_0|H_0)$ to be the $p$-value of a Haar coefficient where $k_0 = 1, 2, \cdots$. Accordingly, a Bi-Haar coefficient can be written as $d_j^{bh} = 2^{-cj}r(X_3-X_4 + \frac{1}{8}(X_1-X_2))$, where $X_1,X_2\sim\cP(\lambda_j)$ and $X_3,X_4\sim\cP(\lambda_j/2)$ are all independent. We note $p_{BH} := \Pr(d^{bh}_j \geq 2^{-cj}k_0|H_0)$ to be the $p$-value of a Bi-Haar coefficient at the same critical threshold as for $p_H$. These definitions can be extended to higher dimensions ($q>1$) straightforwardly.
 
For high-intensity settings or for large scales, $d_j^h$ and $d_j^{bh}$ will be asymptotically normal with the same asymptotic variances $\sigma_h^2 = \sigma_{bh}^2 = 2^{qj(1-2c)}\lambda$ due to the normalized filter banks. Thereby, they will have asymptotically equivalent tail probabilities, i.e., $p_{BH} \approx p_H$.

For low intensity settings ($\lambda \ll 1$) and small scales, the following proposition (proof in Appendix \ref{sec:proof:prop:HaarBiHaarPValues}) shows for 1D signals that $p_{BH}$ is essentially upper-bounded by $p_H$ under $H_0$. The bounds for multidimensional data ($q>1$) are also studied in Appendix~\ref{sec:proof:prop:HaarBiHaarPValues}.
\begin{prop}
\label{prop:HaarBiHaarPValues}
We have the following upper-bound for 1D signals
\begin{equation}
\label{eq:pBHBound}
p_{BH} \leq p_H + A(\lambda_j)(1-2p_H)
\end{equation}
where \[ A(\lambda_j) = \frac{1}{2} \left[1-e^{-2\lambda_j}\left(I_0(2\lambda_j) + 2\sum_{m=1}^8 I_m(2\lambda_j) \right)\right] \]
As $\lambda \to 0+$, $A(\lambda_j) = \frac{2^{9j-7}}{2835}\lambda^9 + o(\lambda^9)$.
\end{prop}
This theoretical bound is clearly confirmed by the numerical simulations shown in Table~\ref{tab:pvalues}. Here we show the results for $\lambda_j\in [10^{-1}, 10^2]$ and different critical thresholds $k_0$ at the tails of the distributions. We indeed observe that $p_{BH}$ is always strictly smaller than $p_H$.
\begin{center}
{\small
\begin{threeparttable}[htbp]
\caption{$p_{H}$ and $p_{BH}$}
\begin{tabular*}{\textwidth}{@{\extracolsep{\fill}}p{0.07\textwidth}|p{0.12\textwidth}|p{0.07\textwidth}|p{0.12\textwidth}|p{0.09\textwidth}|p{0.12\textwidth}|p{0.09\textwidth}|p{0.12\textwidth}}
\hline
\multicolumn{2}{c|}{$\lambda_j = 10^{-1}$} & \multicolumn{2}{c|}{$\lambda_j = 10^{0}$} & \multicolumn{2}{c|}{$\lambda_j = 10^{1}$} & \multicolumn{2}{c}{$\lambda_j = 10^{2}$}\\
\hline
$k_0 = 2$ & ($1.15\times 10^{-3}$, $1.17\times 10^{-4}$) & $k_0 = 4$ & ($1.12\times 10^{-3}$, $4.57\times 10^{-4}$) & $k_0 = 9$ & ($3.97\times 10^{-3}$, $2.48\times 10^{-3}$) & $k_0 = 20$ & ($2.56\times 10^{-2}$, $2.28\times 10^{-2}$) \\
$k_0 = 3$ & ($1.91\times 10^{-5}$, $1.87\times 10^{-6}$) & $k_0 = 5$ & ($1.09\times 10^{-4}$, $4.34\times 10^{-5}$) & $k_0 = 12$ & ($2.12\times 10^{-4}$, $1.26\times 10^{-4}$) & $k_0 = 30$ & ($1.62\times 10^{-3}$, $1.39\times 10^{-3}$) \\
$k_0 = 4$ & ($2.38\times 10^{-7}$, $2.28\times 10^{-8}$) & $k_0 = 6$ & ($8.90\times 10^{-6}$, $3.49\times 10^{-6}$) & $k_0 = 15$ & ($6.60\times 10^{-6}$, $3.78\times 10^{-6}$) & $k_0 = 40$ & ($4.22\times 10^{-5}$, $3.52\times 10^{-5}$) \\
\hline
\end{tabular*}
\vspace{3pt}
\parbox{\textwidth}{Every parenthesis shows $(p_H, p_{BH})$ for 1D signals, where we always observe that $p_{BH} < p_H$.}
\label{tab:pvalues}
\end{threeparttable}
}
\end{center}

\subsection{Thresholds controlling FPR\label{subsec:thresh}}
For individual tests controlling FPR, the HTs can be implemented by thresholding operators. In other words, one can find $\tilde{t}_j$ such that $\Pr(|d^{bh}_j|\geq \tilde{t}_j|H_0) \leq \alpha$ where $\alpha$ represents the controlled FPR. Now consider the Haar case and suppose that we have derived the Haar threshold $t_j$ under the controlled FPR. Then, by setting $\tilde{t}_j := 2^{-cjq}\lceil 2^{cjq}t_j\rceil$ the results in Section~\ref{subsec:bihaar} allow us to conclude that the FPR for a Bi-Haar test will always be upper-bounded by $\alpha$. We point out that to simplify the presentation, $t_j$ and $\tilde{t}_j$ are supposed to be scale-dependent only, but scale \textit{and} location-dependent thresholds can be derived using the same procedure presented below.
 
\subsubsection{CLTB threshold \cite{Kolaczyk1997,Kolaczyk1999a,Kolaczyk2000,Charles2003}\label{subsubsec:cltb}}
The Haar coefficient for $q$D data can be written as $d^h_{j} = 2^{-cjq}(X_1-X_2)$ where $X_1,X_2\sim\cP(\lambda_j/2)$ are independent. It follows from (\ref{eq:haartail}) that:
\begin{eqnarray}
  \label{eq:approxcenchi2} 
  \Pr(d_j^h \geq t_j|H_0) = \Pr\left(\chi_{(2m_j)}^2(\lambda_j) < \lambda_j\right) 
  &\approx& \Pr(\gamma\chi_{(f)}^2 <\lambda_j)\\
  \label{eq:approxclt}
  &\approx& \Pr\left(Z > \frac{f-\lambda_j/\gamma}{\sqrt{2f}}\right)
\end{eqnarray}
where $m_j = 2^{cjq}t_j$, $\gamma = (2m_j+2\lambda_j)/(2m_j+\lambda_j)$,
$f=(2m_j+\lambda_j)^2/(2m_j+2\lambda_j)$, $\chi_{(v)}^2$ is a central
chi-square variable and $Z\sim \mathcal{N}(0, 1)$.
Here, two stages of approximation are used: 1) the non-central chi-square distribution is first approximated by a central one (\ref{eq:approxcenchi2}) \cite{Patnaik1949}; 2) the central chi-square variable is then approximated by a normal one (\ref{eq:approxclt}) using the central
limit theorem (CLT). $t_j$ is thus called the CLT-based (CLTB) threshold.
Consequently, it remains to solve the equation $(\ref{eq:approxclt})=\alpha/2$, and the solution is given by:
\begin{equation}
\label{eq:cltbth}
  t_j = 2^{-cjq-1}\left(z_{\alpha/2}^2+\sqrt{z_{\alpha/2}^4+4\cdot\lambda_j z_{\alpha/2}^2}\right)
\end{equation}
where $z_{\alpha/2} = \Phi^{-1}(1-\alpha/2)$, and $\Phi$ is the standard normal \textit{cdf}. Universal threshold can also be obtained by setting $z_{\alpha/2} = \sqrt{2\ln N_j}$ in (\ref{eq:cltbth}) where $N_j$ is the total number of coefficients in one band at scale $j$. 

\subsubsection{FAB threshold\label{subsubsec:fab}}
An improvement of CLTB threshold can be achieved by replacing (\ref{eq:approxclt}) with an approximation of faster convergence, e.g., the following one proposed by Fisher \cite{Fisher1950}:
\begin{equation}
  \label{eq:approxfisher}
  \sqrt{2\chi_{(f)}^2} \to \mathcal{N}(\sqrt{2f-1}, 1), \quad f \to \infty
\end{equation}
Therefore, (\ref{eq:approxclt}) is changed to:
\begin{eqnarray}
  \Pr\left(\gamma\chi_{(f)}^2 < \lambda_{j}\right) &\approx& \Pr \left(Z > \sqrt{2f-1} - \sqrt{\frac{2\lambda_{j}}{\gamma}}\right)
\end{eqnarray}
Let us denote:
\begin{equation}
  \label{mkolag} 
  G(m_{j}) := \sqrt{2f-1} - \sqrt{\frac{2\lambda_{j}}{\gamma}} 
  = \sqrt{\frac{(2m_{j}+\lambda_{j})^2}{m_{j}+\lambda_{j}}-1} - \sqrt{\frac{\lambda_{j}(2m_{j}+\lambda_{j})}{m_{j}+\lambda_{j}}}
\end{equation}
It remains to solve $G(m_{j}) = z_{\alpha/2}$, which leads to a quartic equation in $m_{j}$:
\begin{eqnarray}
  \nonumber
&16 m_{j}^4 + \left[ 16\lambda_{j} - 8(z_{\alpha/2}^2+1) \right] m_{j}^3 + \left[(z_{\alpha/2}^2+1)^2 - (20z_{\alpha/2}^2 +12)\lambda_{j} + 4\lambda_{j}^2\right]m_{j}^2 & \\
\label{eq:poly4}
& +\left[ 2(z_{\alpha/2}^2+1)^2\lambda_{j} - 16z_{\alpha/2}^2\lambda_{j}^2-4\lambda_{j}^2\right]m_{j} + 
   (z_{\alpha/2}^2+1)^2\lambda_{j}^2 - 4z_{\alpha/2}^2\lambda_{j}^3 = 0&
\end{eqnarray}
The final Fisher-approximation-based (FAB) threshold $t_{j}$ is obtained from $m_j^*$, the solution of (\ref{eq:poly4}). Owing to the following results, we do not need to write out the explicit expression of $m_j^*$, which could be rather complex: 
\begin{prop}
\label{prop:feassol}
The feasible condition for $m_j$ is given by (\ref{eq:feascond}), and the feasible solution $m_j^*$ exists and is unique.
\begin{equation}
\label{eq:feascond}
m_{j} \geq \frac{1}{8} \left[ z_{\alpha/2}^2 - 2\lambda_{j} + 1 + 
      	\left(z_{\alpha/2}^4 + (12\lambda_{j}+2) z_{\alpha/2}^2 + 4\lambda_{j}^2 + 12\lambda_{j} + 1\right)^{1/2}\right]
\end{equation}  
\end{prop}
Proposition~\ref{prop:feassol} implies that we can use any numerical  
quartic-equation solver, e.g.~Hacke's method \cite{Hacke1941}, to find the four solutions of (\ref{eq:poly4}). One and only one of the solutions will satisfy (\ref{eq:feascond}), which is $m_j^*$. The universal threshold can also be derived in the same way as in the CLTB case. 

\subsection{Summary of the denoising controlling FPR\label{subsec:algo}}
Note that the thresholds $\tilde{t}_j$ depend on the background rate at scale $j$ (i.e. $\lambda_{j}$). Without any prior knowledge, it can be estimated by the values of the approximation coefficients at scale $j+1$ (i.e. $a_{j+1}$). Here, the wavelet denoising should be carried out in a coarse-to-fine manner, outlined as follows:
\begin{algorithm}[H]
\caption{Poisson noise removal by HTs in the Bi-Haar domain}
\begin{algorithmic}[1]
\STATE Bi-Haar transform of $\fv$ up to $j=J$ to obtain $a_{J}$ and $d^{bh}_{j}$~($1 \leq j \leq J$)
\FOR{$j=J$ down to $1$}
		\STATE $\hat{\lambda}_{j} = 2^{jq}\lambda$ if $\lambda$ is known; otherwise $\hat{\lambda}_{j} = \max(2^{cjq}a_{j}, 0)$
		\STATE Testing $d^{bh}_{j}$ by applying thresholds $\tilde{t}_j$ for a prefixed $\mbox{FPR} = \alpha$
		\STATE Reconstruct $a_{j-1}$ by inverse Bi-Haar transform
\ENDFOR
\STATE Positivity projection: $\hat{\Lambda} = \max(a_{0}, 0)$
\end{algorithmic}
\end{algorithm}

\section{Results\label{sec:res}}
\subsection{Haar vs. Bi-Haar denoising for regular intensities\label{subsec:hbh}}
\begin{figure}[htbp]
\centering
\mbox{
      \subfigure[]{\epsfig{file=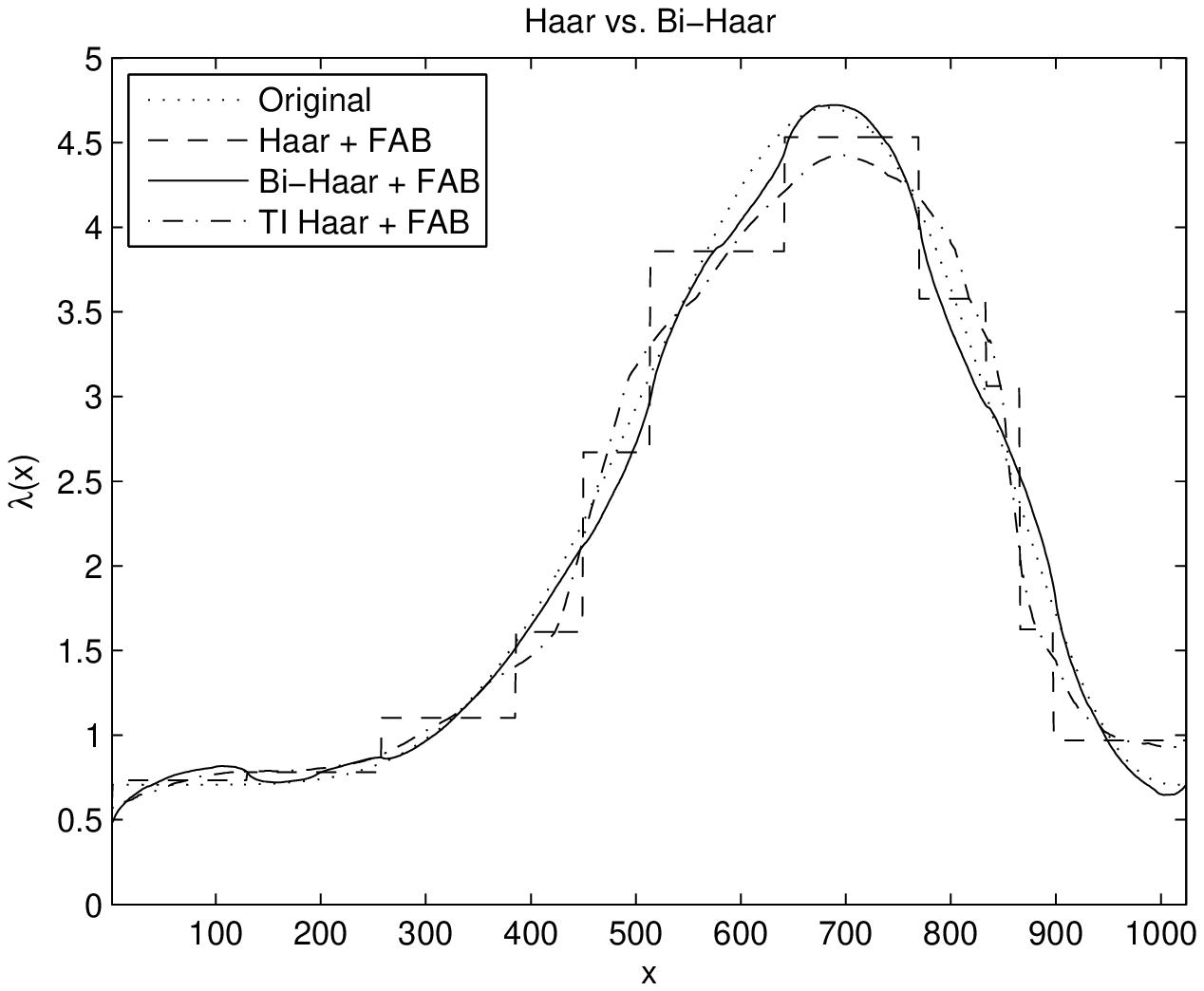, width = 0.45\textwidth, clip}}
      \subfigure[]{\epsfig{file=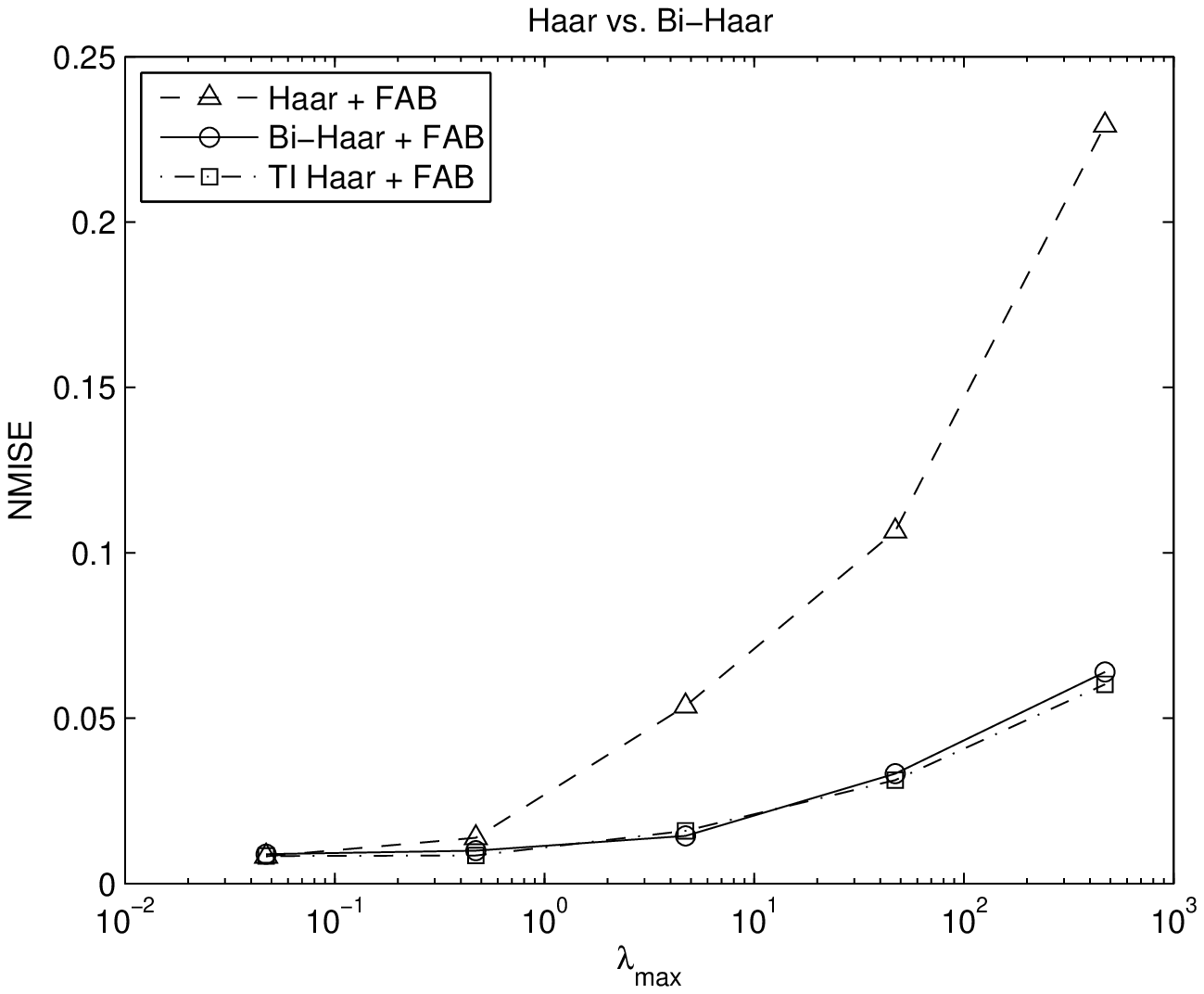, width=0.45\textwidth, clip}}
}
  \caption[Denoising the ``Smooth'' function]{Denoising the ``Smooth'' function (length = 1024). Estimates from Haar, Bi-Haar and TI Haar (undecimated) are compared. $\alpha = 10^{-3}$ and $J = 7$. (a) denoising results; (b) NMISEs.}
  \label{fig:HaarvsBiHaar}
\end{figure}

To compare Haar and Bi-Haar denoising for regular intensities, 
we generate noisy signals from the ``Smooth'' function \cite{Besbeas2004} (see Fig.\ref{fig:HaarvsBiHaar}(a)) and measure the Normalized Mean Integrated Square Error (NMISE) per bin from the denoised signals. The NMISE is defined as: $\mbox{NMISE} := \bE[(\sum_{i=1}^N (\hat{\lambda}_i - \lambda_i)^2/\lambda_i)/N]$, where $(\hat{\lambda}_i)_i$ is the intensity estimate. Note that the denominator $\lambda_i$ plays the role of variance stabilization in the error measure.

Fig.\ref{fig:HaarvsBiHaar}(a) shows the denoising examples given by Haar, Bi-Haar and TI Haar estimations, where FAB thresholds are applied to control a FPR $\alpha=10^{-3}$. The original intensity function is scaled to cover a wide range of intensities, and Fig.\ref{fig:HaarvsBiHaar}(b) compares the NMISEs (measured from $100$ replications) of the three estimators as functions of the underlying peak intensity. 

It can be seen that the Bi-Haar estimate is much more regular than the Haar one, and is even almost as good as TI Haar at every intensity level under the NMISE criterion. This surprising performance is gained with the same complexity as in the Haar denoising, i.e., $O(N)$ only, as opposed to $O(N\log N)$ in the TI Haar case. 

\subsection{Source-flux estimation in astronomical hyperspectral data\label{subsec:astro}}
We apply our method to source-flux estimation in astronomical hyperspectral images.
A hyperspectral image $\fv(x,y,\nu)$ is a ``2D+1D'' volume, where $x$ and
$y$ define the spatial coordinates and $\nu$ indexes the spectral band. Each bin records the detected number of photons. As the three axes of
our data have different physical meanings, we are motivated to apply a ``2D+1D'' wavelet transform instead of using the classical 3D transform. That is, we first carry out a complete 2D wavelet
transform for spatial planes, and then a 1D transform along
the spectral direction. We use $j_{xy}$ and $j_{\nu}$ to denote the $j$-th spatial scale and the $j$-th spectral scale, respectively. Hyperspectral data in practice can be very large, implying that fast denoising is only possible with decimated transforms (the execution time of the example below on a P4 2.8GHz PC is $13$s for our Bi-Haar denoising, i.e., more than $50$ times faster than the TI Haar denoising ($665$s)), not to mention the memory space required by the redundant TI transform.

\begin{figure}[htbp]
  \centering
    \mbox{
      \subfigure[]{\epsfig{file=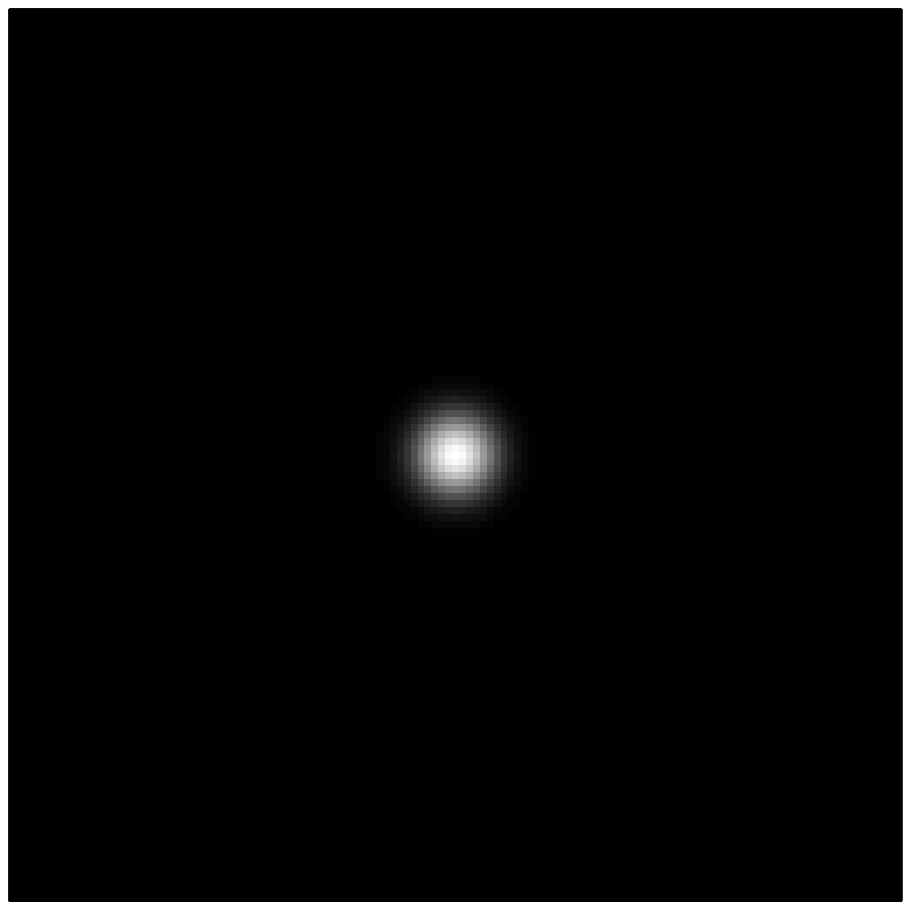, width = .32\textwidth, clip}}
      \subfigure[]{\epsfig{file=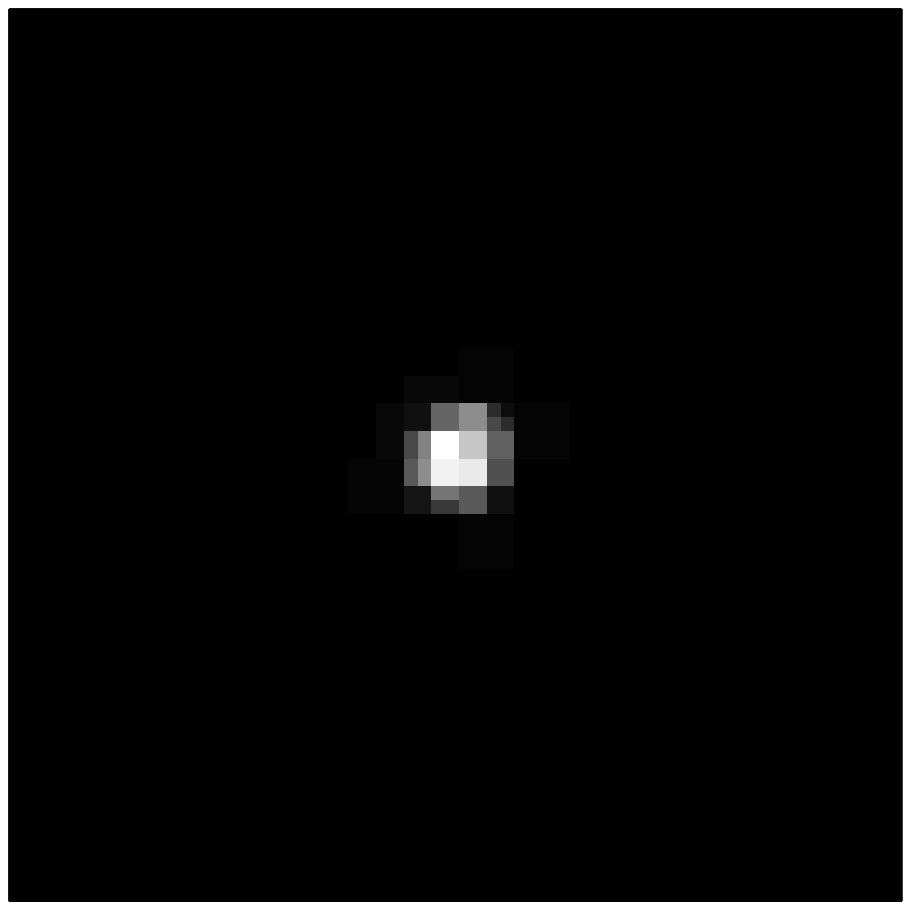, width = .32\textwidth, clip}} 
      \subfigure[]{\epsfig{file=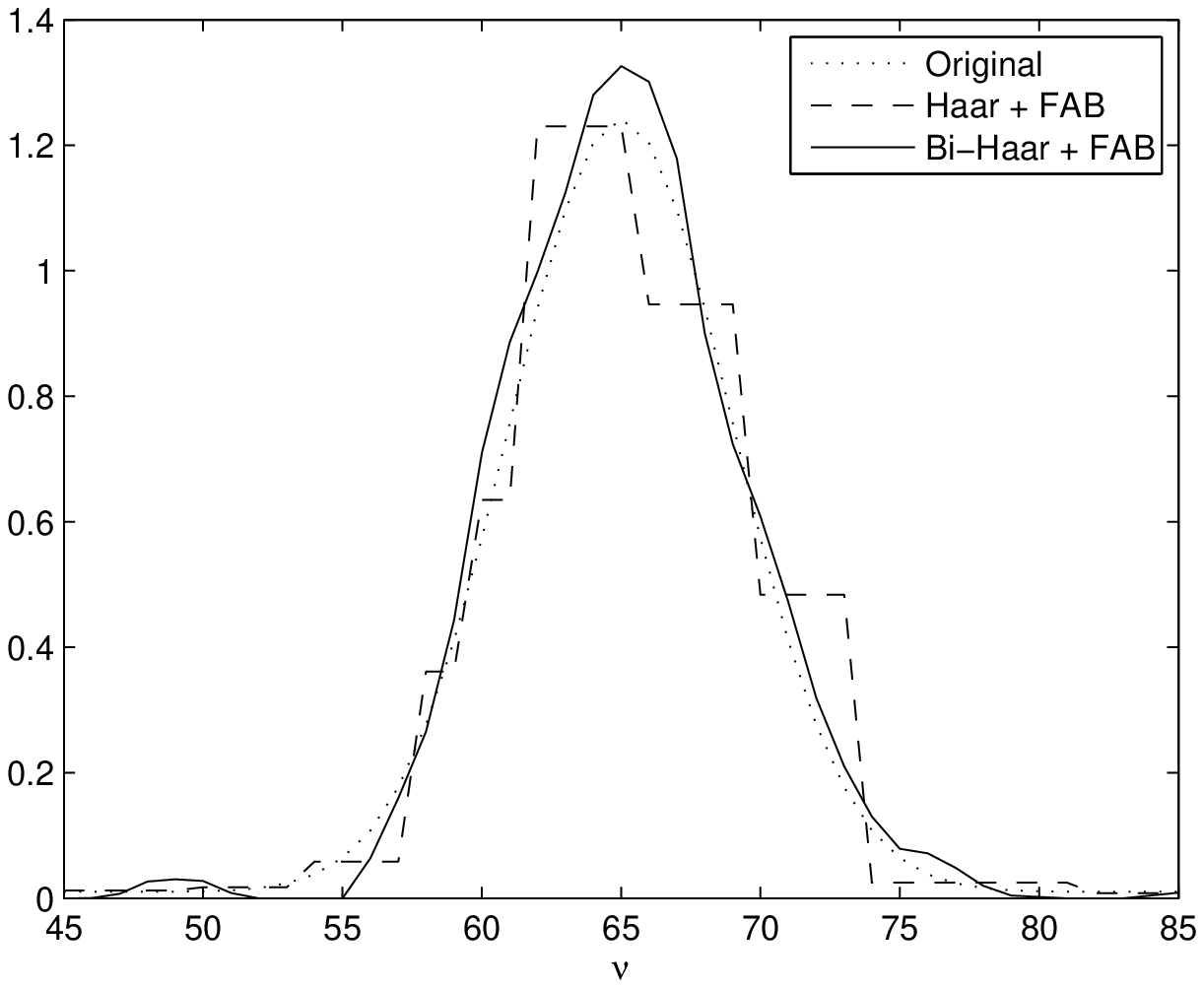, width = .32\textwidth, clip}} 
      }
    \mbox{
      \subfigure[]{\epsfig{file=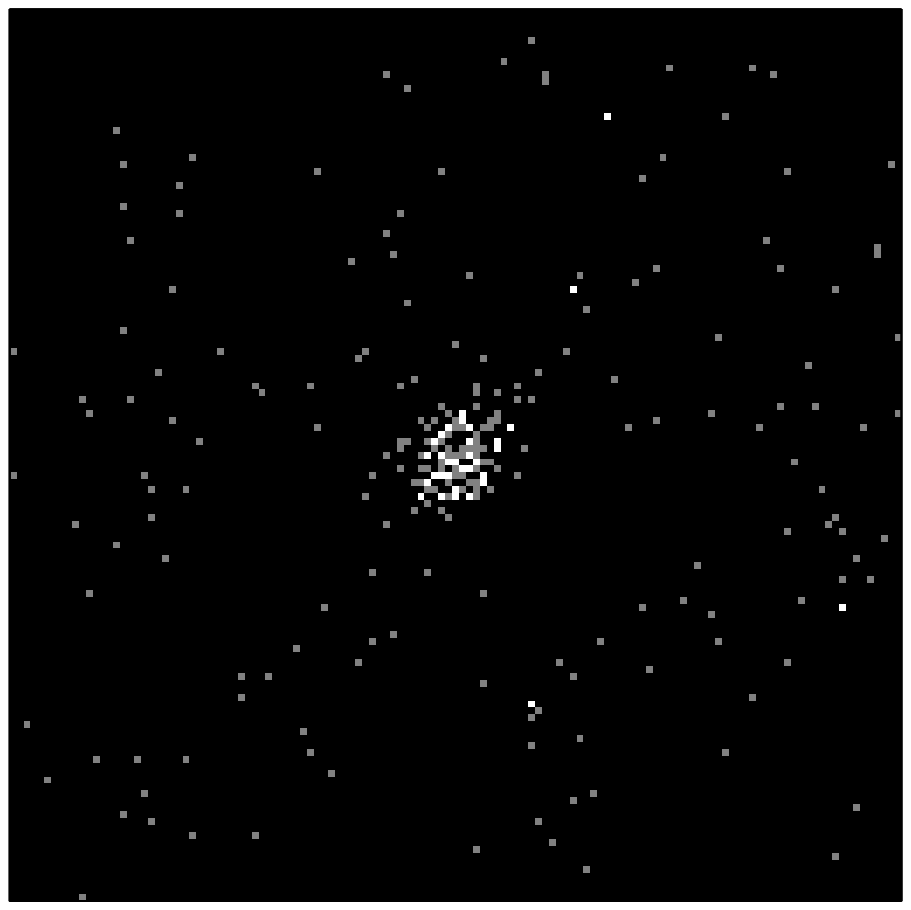, width = .32\textwidth, clip}}
      \subfigure[]{\epsfig{file=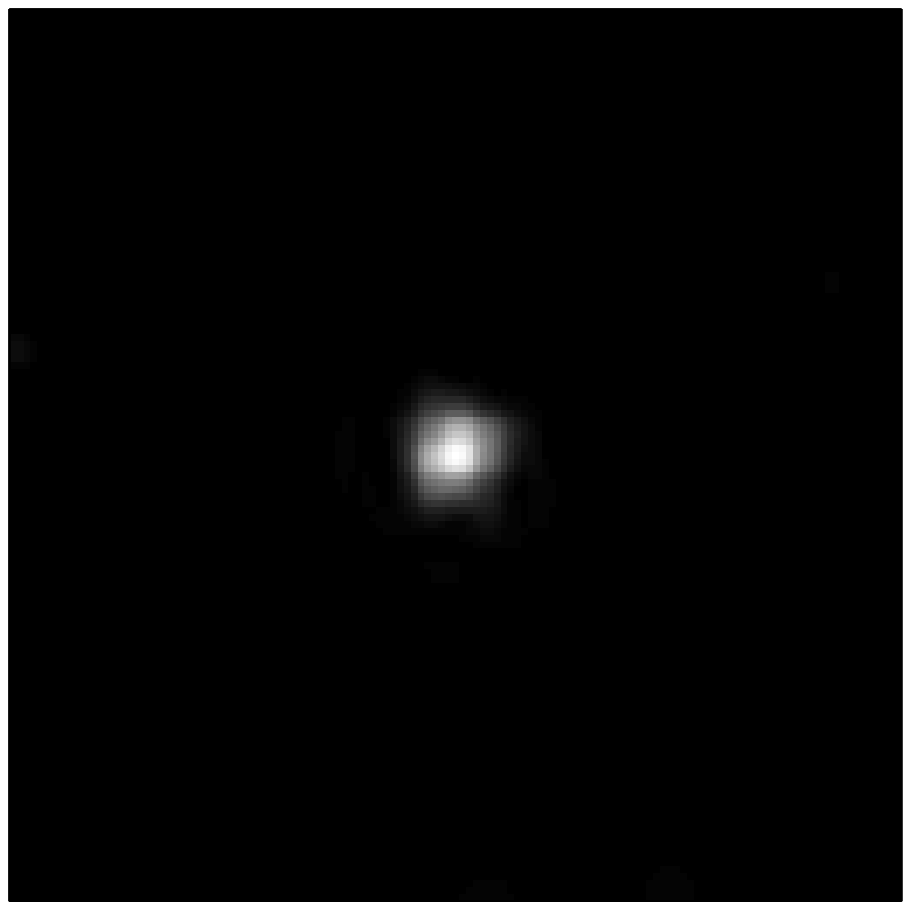, width = .32\textwidth, clip}} 
      \subfigure[]{\epsfig{file=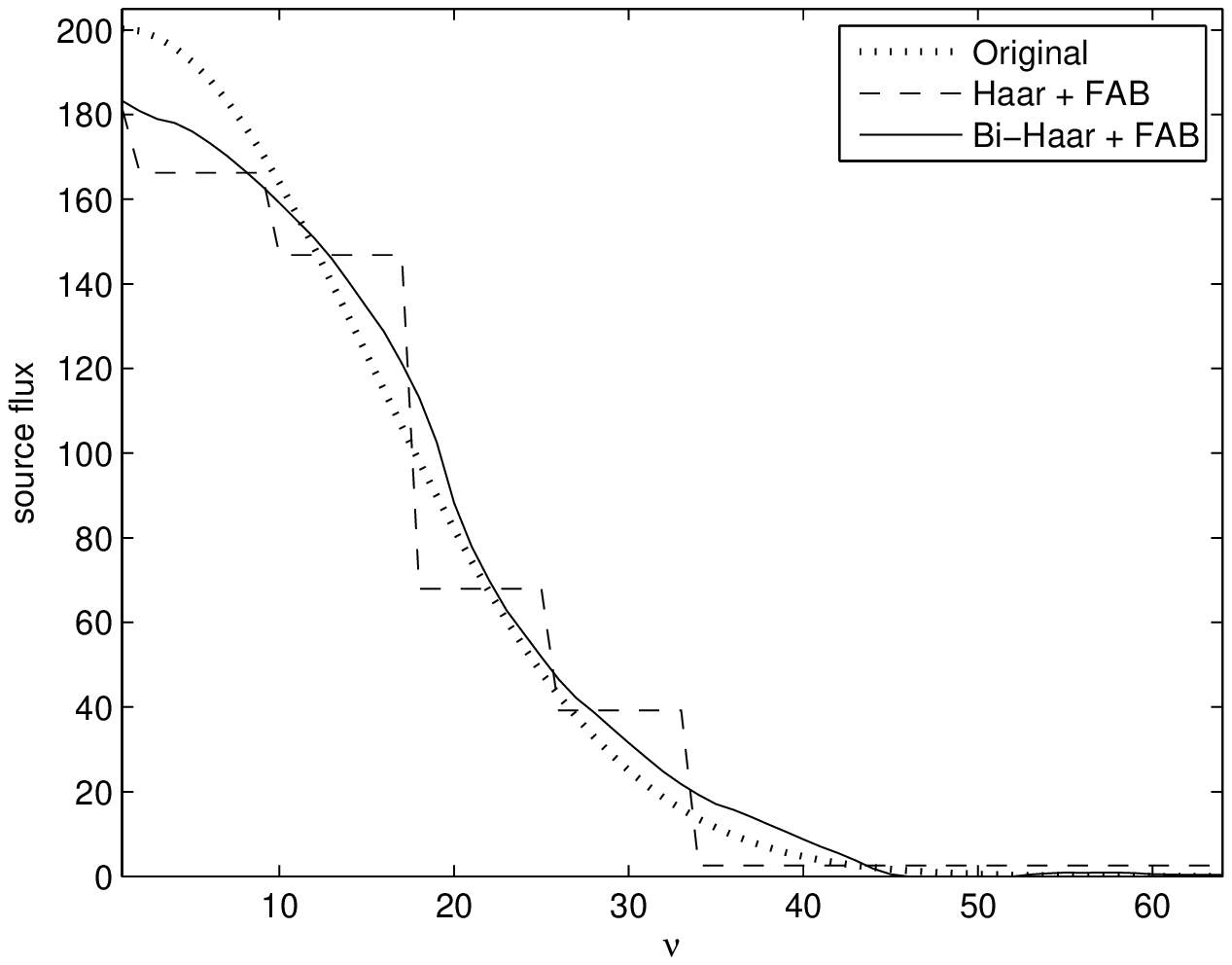, width = .32\textwidth, clip}} 
      }
  \caption{Source-flux estimation in a hyperspectral image (size: $129\times 129\times 64$). $A_\nu \in [10^{-4}, 2]$; $J_{xy}=3$, $J_\nu=5$, FAB thresholding with $\alpha = 10^{-5}$. (a) intensities at $\nu=15$; (b) Haar-denoised data ($\nu=15$); (c) estimated source profile at $\nu=15$ (intensity along a line passing through the source center);  (d) Poisson count image; (e) Bi-Haar-denoised data ($\nu=15$); (f) estimated flux (Respectively for Haar and Bi-Haar estimates: $\frac{1}{\sqrt{N}}\|\hat{S}-S\|_{\ell^2} = 14.4$ and $7.7$.)} 
\label{fig:hypersp}
\end{figure}

Our simulated data contain a source having a Gaussian profile. The source amplitude $A_\nu$ decreases from $2$ to $10^{-4}$ as $\nu$ increases. One example band is shown in Fig.\ref{fig:hypersp}(a). The observed counts at that band are depicted in Fig.\ref{fig:hypersp}(d). The denoising results using Haar and Bi-Haar transforms are respectively shown in Fig.\ref{fig:hypersp}(b) and (e), where FAB thresholds are applied. Fig.\ref{fig:hypersp}(c) illustrates the estimation smoothness gained by Bi-Haar by comparing a line profile of the estimated source from different methods. 
In hyperspectral imaging, the source flux $S(\nu)$ is an important quantity, which equals to the integral of the source intensity over its spatial support at band $\nu$.
Fig.\ref{fig:hypersp}(f) compares the flux given by different denoisers. Clearly, the Haar-based approach leads to a piecewise constant estimate, whereas Bi-Haar provides a regular flux which is more accurate: the normalized $\ell^2$-loss for Haar and Bi-Haar flux estimates, i.e.~$\frac{1}{\sqrt{N}}\|\hat{S}-S\|_{\ell^2}$, are $14.4$ and $7.7$ respectively. 

\section{Conclusion}
\label{sec:conc}
In this paper, we proposed to combine the HT framework with the decimated Bi-Haar transform instead of the classical Haar for denoising large datasets of Poisson counts. We showed that the Haar-based individual HTs can be applied to Bi-Haar coefficients to control a prefixed FPR. By doing so, we benefit from the regular Bi-Haar filter bank to gain a smooth estimate with no ``staircase'' artifacts, while always maintaining a low computational complexity. A Fisher-approximation-based threshold implementing HTs is also designed. This approach could be extended in the future to fast deconvolution of Poisson data.
 
\appendix
\section{Proof of Proposition~\ref{prop:HaarBiHaarPValues}\label{sec:proof:prop:HaarBiHaarPValues}}
\begin{pf}
We note that  
\[ p_H = \Pr(d_j^h \geq 2^{-cj} k_0|H_0) = \sum_{k \geq k_0} e^{-\lambda_j} I_k(\lambda_j) = \sum_{k \leq -k_0} e^{-\lambda_j} I_{|k|}(\lambda_j) \]
where $k_0 \geq 1$. The $p$-value of $d_{BH}$ is given by
\begin{eqnarray*}
p_{BH} &=& \Pr(X_1-X_2 + 8(X_3-X_4) \geq \lceil 8k_0/r\rceil|H_0)  \\
&=& \sum_{k\in\bZ} \Pr(X_3-X_4=k|H_0)\sum_{n=\lceil 8k_0/r\rceil}^\infty\Pr(X_1-X_2 = n-8k|H_0)
\end{eqnarray*}
where $X_{1,2}\sim \cP(\lambda_j)$, $X_{3,4}\sim \cP(\lambda_j/2)$, and $(X_i)_i$ are independent. Now we have,
\begin{eqnarray}
\nonumber
p_{BH} &=& \sum_{k\geq k_0} \Pr(X_3-X_4=k|H_0) \cdot \\\nonumber
&&\sum_{n=\lceil 8k_0/r\rceil}^\infty \Pr(X_1-X_2 = n-8k|H_0)+\Pr(X_1-X_2 = n+8k|H_0) \\
\label{eq:pBHpart2}
&+& \sum_{|k|< k_0} \Pr(X_3-X_4=k|H_0)\sum_{n=\lceil 8k_0/r\rceil}^\infty \Pr(X_1-X_2 = n-8k|H_0)\\
\nonumber
&\leq& p_H + \sum_{|k|<k_0} e^{-\lambda_j} I_{|k|}(\lambda_j) \underbrace{ \sum_{n=\lceil 8k_0/r\rceil}^\infty e^{-2\lambda_j} I_{|n-8k|}(2\lambda_j)}_{T}
\end{eqnarray}
To bound $T$, we use the identity \cite{Abramowitz1970}
$e^x = I_0(x) + 2\sum_{n=1}^\infty I_n(x)$. 
As $r < 1$, we have
\begin{eqnarray*}
T &\leq& e^{-2\lambda_j}\sum_{n\geq 9} I_n(2\lambda_j) = \frac{1}{2} \left[1-e^{-2\lambda_j}\left(I_0(2\lambda_j) + 2\sum_{m=1}^8 I_m(2\lambda_j) \right)\right] =: A(\lambda_j)
\end{eqnarray*}
Thus, $p_{BH} \leq p_H + A(\lambda_j)(1-2p_H)$. 
As $\lambda \to 0+$, we have that $A(\lambda_j) = \frac{2^{9j-7}}{2835}\lambda^9 + o(\lambda^9)$. \qed
\end{pf}

The same arguments can be carried out to bound $p_{BH}$ in multi-dimensional cases. As an example, let us consider 2D data. A 2D wavelet transform will produce bands of horizontal, vertical and diagonal detail coefficients, i.e., $d_{j;H}$, $d_{j;V}$, and $d_{j;D}$. For horizontal and vertical coefficients, we have that $p_{BH} \leq p_H + A(\lambda_j)(1-2p_H)$, where $\lambda_j := 4^j\lambda$. For diagonal coefficients, it can be shown that $p_{BH} \leq p_H + B(\lambda_j) (1-2p_H)$, where
\begin{equation}
\label{eq:Bj2D}
B(\lambda_j) := 
\frac{1}{2} \left[1 - \sum_{n=-64}^{64} \sum_{k=-\infty}^\infty e^{-8\lambda_j}I_{|k|}(4\lambda_j)I_{|n-8k|}(4\lambda_j)\right]
\end{equation}
To see the behavior of $B(\lambda_j)$ as the intensity becomes small, we note $B(\lambda_j) = B_K(\lambda_j) + \epsilon_K$. Here, $B_K$ is given by (\ref{eq:Bj2D}) with $k$ ranging from $-K$ to $K$, and $\epsilon_K$ is the residual which can be made arbitrary small as $K$ increases. Then, we have for all $K\geq 8$ that $B_K(\lambda_j) = \frac{8}{567} \lambda_j^9 + o(\lambda_j^9)$. Clearly, this procedure can be continued for higher dimensional cases ($q>2$). 


\section{Proof of Proposition~\ref{prop:feassol}\label{sec:prf:prop:feassol}}
\begin{pf}
The facts that $G(m_{j}) = z_{\alpha/2}$, $z_{\alpha/2} >0$, 
$2f - 1 \geq 0$, $m_{j} > 0$ and $\lambda_{j} \geq 0$ show (\ref{eq:feascond}). 

Next, when the equality in (\ref{eq:feascond}) holds, we have:
\[ G(m_{j}) = \sqrt{z_{\alpha/2}^2 + \frac{\lambda_j}{2m_{j}}(z_{\alpha/2}^2+1)} - \sqrt{\frac{\lambda_j}{2m_{j}}(z_{\alpha/2}^2+1)} \leq  z_{\alpha/2} \]
The existence and uniqueness of the feasible solution follow from the fact that $G$ is a strictly increasing function under (\ref{eq:feascond}), and that $G(m_{j}) \to +\infty$ as $m_{j} \to +\infty$. \qed
\end{pf}

\end{document}